\documentclass[11pt]{amsart}

\usepackage[T1]{fontenc}
\usepackage[utf8]{inputenc}
\usepackage{amsmath,amssymb,amsthm,mathtools,xcolor}
\usepackage{enumitem}
\usepackage{hyperref}
\hypersetup{hidelinks}
\usepackage[margin=1.15in]{geometry}
\allowdisplaybreaks

\theoremstyle{plain}
\newtheorem{theorem}{Theorem}[section]
\newtheorem{proposition}[theorem]{Proposition}

\newtheorem{corollary}[theorem]{Corollary}

\theoremstyle{definition}

\newtheorem*{theorem*}{Theorem}

\theoremstyle{remark}
\newtheorem{remark}[theorem]{Remark}

\numberwithin{equation}{section}

\begin{document}

\title[Harmonic-Einstein maps and critical metrics for eigenvalues]{Harmonic-Einstein maps and critical metrics for eigenvalues of the conformal Laplacian}

\author{Mehdi Eddaoudi}
\address{Department of Mathematics and Data Science, Vrije Universiteit Brussel.
}
\email{mehdi.eddaoudi@vub.be}

\begin{abstract}
We prove that Clarke critical metrics for the ratio of the first two eigenvalues of the conformal Laplacian on a closed manifold with positive Yamabe constant correspond to harmonic--Einstein maps into a sphere. We discuss some rigidity consequences, as well as an extension to the higher-order GJMS operators.
\end{abstract}

\maketitle
\tableofcontents

\section{Introduction}

\subsection{Eigenvalue ratio of the conformal Laplacian}

Let $(M^n,g)$ be a smooth closed connected Riemannian manifold
of dimension $n\geq3$, and let $L_g$ be the conformal Laplacian defined as
\[
L_g=\Delta_g+c_nS_g,
\qquad
c_n=\frac{n-2}{4(n-1)},
\]
where $\Delta_g=-\operatorname{div}_g\nabla$ is the nonnegative
Laplace--Beltrami operator and $S_g$ denotes the scalar curvature.
It is well known that the spectrum of $L_g$ consists of a discrete sequence of real eigenvalues
\[
\lambda_1(M,L_g)<\lambda_2(M,L_g)\leq\cdots\longrightarrow+\infty,
\]
where each eigenvalue is repeated according to its multiplicity.
In particular, the first eigenvalue is always simple, and a corresponding eigenfunction can be chosen strictly positive.

Since the eigenvalues of $L_g$ are not scale-invariant,
\[
\lambda_k(M,L_{c^2g})=\frac{1}{c^2}\,\lambda_k(M,L_g), \qquad c>0,
\]
any variational approach must rescale them by
a suitable geometric quantity. The most common
choice is the volume $\operatorname{Vol}(M,g)$. Within a conformal class, minimizing
$\lambda_1(M,L_g)\operatorname{Vol}(M,g)^{2/n}$ is the spectral
formulation of the Yamabe problem
\cite{Yamabe1960Deformation,Trudinger1968ConformalDeformation,Aubin1976Yamabe,Schoen1984Yamabe,LeeParker1987Yamabe}.
For higher eigenvalues, Ammann and Humbert \cite{AmmannHumbert2006SecondYamabe}
proved under some assumptions the existence of a generalized metric minimizing the second volume-normalized eigenvalue
\(\lambda_2(M,L_g)\operatorname{Vol}(M,g)^{2/n}\)
within a conformal class of positive Yamabe type. El Sayed \cite{ElSayed2014} then obtained further results on nodal solutions of
Yamabe-type equations, while Gursky and Pérez-Ayala
\cite{GurskyPerezAyala2022SecondEigenvalue} treated conformal classes for
which $L_g$ has at least two negative eigenvalues, proving that extremal
metrics for the second eigenvalue yield either a nodal solution of the
Yamabe equation or a harmonic map into a sphere. More recently, Humbert, Petrides and Premoselli
\cite{HumbertPetridesPremoselli2025} investigated extremal
metrics for eigenvalues of any index of the subcritical GJMS operators
within a fixed conformal class. They carried out a detailed analysis of
these extremals, proving existence results and deriving the
Euler--Lagrange equations satisfied by local extremals.

On the other hand, the maximization problem of volume-normalized eigenvalues behaves differently. Ammann and Jammes
\cite{AmmannJammes2011Supremum} proved that the first positive
volume-normalized eigenvalue is unbounded from above, even within a fixed
conformal class. This motivated the search for other normalizations, such
as the one introduced by Sire and Xu \cite{SireXu}, for which the supremum
over a conformal class is finite. Pérez-Ayala
\cite{PerezAyala2021SireXu} then extended this normalization to the $k$-th
eigenvalue, obtained necessary conditions for the existence of extremal
metrics, and settled the existence problem for $k=1$.

In this paper, we propose yet another normalization: when
$\lambda_1(M,L_g)>0$, we divide by the first eigenvalue itself and consider
the ratio
\[
\mathcal R_2(g)=\frac{\lambda_2(M,L_g)}{\lambda_1(M,L_g)}.
\]
Observe that the ratio $\mathcal R_2$ is invariant under constant rescalings of the
metric, therefore finding its critical points is equivalent to finding critical points
of $\lambda_2(M,L_g)$ under the constraint $\lambda_1(M,L_g)=1$. Ratios of
this type already appear in the inequalities of Payne, P\'olya and
Weinberger \cite{PaynePolyaWeinberger1956Ratio} for the Euclidean
Dirichlet problem, and in the work of Ashbaugh and Benguria
\cite{AshbaughBenguria1991,AshbaughBenguria1992Second,AshbaughBenguria1992SharpBound,ashbaugh1994isoperimetric}.
The metric is allowed to vary arbitrarily, not only within a conformal class.

Let $\mathcal M(M)$ be the space of smooth Riemannian metrics on $M$, and let
\[
\mathcal M^+(M)
=\{g\in\mathcal M(M):\lambda_1(M,L_g)>0\}.
\]
Since the sign of $\lambda_1(M,L_g)$ is a conformal invariant and agrees
with the sign of the Yamabe constant $Y(M,[g])$, $\mathcal M^+(M)$ is exactly
the space of metrics whose conformal class is of positive Yamabe type. Our
aim is to describe the geometric properties of critical metrics of
$\mathcal R_2$ on $\mathcal M^+(M)$, and to deduce geometric and topological
restrictions that this criticality imposes.

The main difficulty in this context is that \(\mathcal R_2\) is not
differentiable in general, because of the possible multiplicity of
\(\lambda_2\). For the volume-normalized Laplace eigenvalues
\[
\overline\lambda_k(M,\Delta_g)
=\lambda_k(M,\Delta_g)\operatorname{Vol}(M,g)^{2/n},
\]
a notion of criticality adapted to this lack of differentiability was
introduced by Nadirashvili \cite{Nadirashvili1996Berger} for surfaces, and
extended by El Soufi and Ilias \cite{ElSoufiIlias2008} to higher dimensions.
A metric \(g\) is critical for \(\overline\lambda_k\) if, for every
analytic deformation \(g_t\) with \(g_0=g\),
\[
\left.\frac{d}{dt}\overline\lambda_k(M,\Delta_{g_t})\right|_{t=0^-}
\cdot
\left.\frac{d}{dt}\overline\lambda_k(M,\Delta_{g_t})\right|_{t=0^+}
\leq0.
\]
The resulting Euler--Lagrange equations relate critical metrics to minimal
isometric immersions into spheres when the whole space of metrics is
varied, and to harmonic maps into spheres when only the conformal class is
varied. These ideas have since developed into a rich literature, ranging from
equivariant spectral optimization problems to the construction of minimal
surfaces
\cite{KarpukhinStern2024MinMax,Petrides2023Combinations,Vinokurov2025,vinokurov2026eigenvalueoptimizationfirstvariationformula,VinokurovPHarmonic2026,PetridesRegularity2026,KarpukhinStern2024,KNPP2021,KNPP2,GNP,ColboisElSoufi2003,JakobsonNadirashviliPolterovich2006,GirouardKarpukhinLagace2021,Petrides2014,KarpukhinKusnerMcGrathStern2024,KarpukhinPetridesStern2025}.

Recently, Petrides and Tewodrose \cite{PetridesTewodrose2024} introduced a
notion of criticality based on Clarke's nonsmooth analysis
\cite{Clarke1990Optimization}, which is more unifying and natural than the
one of Nadirashvili and El Soufi and Ilias. It
applies to locally Lipschitz functionals depending on finitely many
eigenvalues, admits a chain rule, and expresses the Euler--Lagrange
equations as a balance between first-variation tensors. This framework is
perfectly adapted to our ratio \(\mathcal R_2\), and it is the one we use
in this paper.

The key idea of this work is to notice that, once the Euler--Lagrange
identity is split into its trace and trace-free parts, the conformal
covariance of $L_g$ unravels a remarkable geometric structure, which we
now describe. Let $u$ be the positive $L^2$-normalized first eigenfunction of $L_g$ and set
$g_u=u^{4/(n-2)}g$. The trace part gives the conformal part of the
equation, which is carried by a harmonic map into a sphere with respect to
$g_u$, while the trace-free part contains the Einstein--harmonic map
equation relating the pull-back of the round metric to the Ricci tensor of
$g_u$. 
\subsection{Main Result}
This is made precise in the following theorem.

\begin{theorem}\label{Thm:principal}
A metric $g\in\mathcal M^+(M)$ is Clarke critical for
$\mathcal R_2$ if and only if there exist a first eigenfunction
$u>0$ and pairwise $L^2(M,g)$-orthogonal second eigenfunctions $v_1,\ldots,v_m$ such that the map
$\psi:M\to\mathbb S^{m-1}$ with coordinates $\psi_j=v_j/u$
is harmonic with respect to $g_u=u^{4/(n-2)}g$ and satisfies
\begin{equation}\label{eq:Einstein-harmonic}
\psi^*g_{\mathbb S^{m-1}}
=c_n\bigl(\mathcal R_2(g)-1\bigr)\operatorname{Ric}_{g_u},
\end{equation}
where $g_{\mathbb S^{m-1}}$ is the standard metric on the unit sphere.
\end{theorem}

Taking the trace in \eqref{eq:Einstein-harmonic}, we obtain
$|d\psi|^2_{g_u}=c_n(\mathcal R_2(g)-1)S_{g_u}$. Together with the
harmonicity of $\psi$, we obtain the trace reversed form
\begin{equation}\label{eq:sigma-Einstein}
\operatorname{Ric}_{g_u}-\dfrac12 S_{g_u}\,g_u
=\dfrac{1}{c_n(\mathcal R_2(g)-1)}
\Bigl(\psi^*g_{\mathbb S^{m-1}}-\dfrac12|d\psi|^2_{g_u}\,g_u\Bigr).
\end{equation}
The left-hand side is the Einstein tensor of \(g_u\), while
the bracket on the right is the stress--energy tensor of \(\psi\). Therefore, this identifies
\((g_u,\psi)\) as a solution of the Riemannian Einstein--sigma model with zero cosmological constant.
We refer to Anselli \cite{Anselli2021BachEinstein} for its variational
characterization, and Colombo, Mariani and Rigoli 
\cite{ColomboMarianiRigoli2024} for rigidity results under
additional curvature assumptions.
For the classical Lorentzian formulation, see Choquet-Bruhat \cite[Sections~III.6.5 and~III.7.2--III.7.3]{ChoquetBruhat2009}.

The pair $(g_u,\psi)$ can also be seen as a stationary solution of the
Ricci--harmonic map flow with a fixed coupling constant $\alpha>0$,
introduced by M\"uller \cite{Muller2012CoupledFlow},
\[
\frac{\partial}{\partial t}g=-2\operatorname{Ric}_g
+2\alpha\,\psi^*g_{\mathbb S^{m-1}},
\qquad
\frac{\partial}{\partial t}\psi=\tau_g(\psi),
\]
where $\tau_g(\psi)$ denotes the tension field of $\psi$. Here
$\alpha=\bigl(c_n(\mathcal R_2(g)-1)\bigr)^{-1}$.

See also
\cite{Lott2010DimensionalReduction,Williams2015CoupledFlows}
for related developments.

\subsection{Rigidity and topological consequences}
Since \(\mathcal R_2(g)>1\), the tensorial equation yields
\[
\operatorname{Ric}_{g_u}\geq0,
\qquad
\operatorname{rank}\operatorname{Ric}_{g_u}
=\operatorname{rank}d\psi.
\]
Combining this relation with the Bochner formula for \(\psi\),
we show in Corollaries~\ref{cor:spherical-rigidity-intro}
and~\ref{cor:rank-r-intro} that
\(r=\max_M\operatorname{rank}\operatorname{Ric}_{g_u}\geq2\)
and we establish the sharp bound
\[
\mathcal R_2(g)
\geq1+\frac{4(n-1)}{(n-2)(r-1)}
\geq\frac{n+2}{n-2}.
\]
Equality in the first estimate forces the universal
Riemannian cover of $(M,g_u)$ to split as a product
of a round $r$-sphere and a Ricci-flat manifold.
Equality in the dimension-dependent estimate forces
$(M,g_u)$ to be a round sphere.
The same rank information also bounds the multiplicity
of the second eigenvalue:
\[
\dim E_{\lambda_2}(g)\geq r+1\geq3.
\]

The nonnegativity of $\operatorname{Ric}_{g_u}$ also implies
that harmonic one-forms are parallel.
By the Cheeger--Gromoll result
\cite{CheegerGromoll1971}, the fundamental group is virtually
abelian. If $q$ denotes the rank of a finite-index free
abelian subgroup, we obtain
\[
b_1(M)\leq q\leq n-r\leq n-2,
\qquad
\mathcal R_2(g)
\geq1+\frac{4(n-1)}{(n-2)(n-q-1)}.
\]
This last estimate is sharp, with equality realized
by products of round spheres and suitably scaled flat tori.

These restrictions imply obstructions to the existence
of critical metrics. For instance, we prove in Corollary \ref{cor:further-topological-obstructions-intro} that no closed aspherical manifold admits
such a critical metric, and if the Euler characteristic $\chi(M)\neq0$, its fundamental
group must be finite.
In particular, we obtain in Corollary \ref{cor:connected-sum-dichotomy-intro} that the connected sums
$\#_k(\mathbb S^{n-1}\times\mathbb S^1)$ with $k\geq2$ admit metrics of positive scalar curvature but no critical
metric for $\mathcal R_2$.

\subsection{Extension to the GJMS operators}

Since the conformal Laplacian is the first member of the GJMS family,
it is natural to ask whether Theorem~\ref{Thm:principal} extends
to higher orders.
The GJMS operators $P_g^s$ of order $2s<n$, introduced by
Graham, Jenne, Mason and Sparling
\cite{GrahamJenneMasonSparling1992}, are classical operators
in conformal geometry. Together with their associated
$Q$-curvatures, they extend the conformal Laplacian and scalar
curvature to higher orders and provide a natural framework
for higher-order analogues of the Yamabe problem.
They include the conformal Laplacian for $s=1$ and
the Paneitz operator for $s=2$ \cite{Paneitz2008,Branson1985}.
We refer to Humbert, Petrides and Premoselli \cite{HumbertPetridesPremoselli2025} for recent results on volume-normalized eigenvalues in a fixed conformal class.

We consider the ratio
\[
\mathcal R_{2s}(g)
=\frac{\lambda_2(M,P_g^s)}{\lambda_1(M,P_g^s)}.
\]
The previous arguments used for the conformal Laplacian rely essentially on its conformal covariance, a first positive and simple eigenvalue, alongside with a positive first eigenfunction. Since a higher-order
operator does not in general satisfy a maximum principle, we therefore shall assume explicitly that $P_g^s$ is coercive,
that its first eigenvalue is simple, and that a corresponding
first eigenfunction is strictly positive.

To state the result, we briefly introduce the notation
and refer to Section~4 for the precise definitions
and variational formulas.
The associated $Q^{(s)}$-curvature and reduced GJMS operator
are given by
\[
P_g^s(1)=\frac{n-2s}{2}Q_g^{(s)},
\qquad
L_g^{(s),0}=P_g^s-P_g^s(1),
\]
where $P_g^s(1)$ acts by multiplication.
For a sphere-valued map $\psi$, let $H_g^{(2s)}(\psi)$
be the metric stress tensor of the reduced energy
\[
\mathcal E^{(2s),0}(\psi,g)
=\sum_{\alpha=1}^m\int_M
\psi_\alpha L_g^{(s),0}\psi_\alpha\,dv_g,
\]
with $\psi$ kept fixed under metric variations.
The metric variation tensor of the total $Q^{(s)}$-curvature
is denoted by $E_g^{Q^{(s)}}$ and satisfies
\[
E_g^{Q^{(s)}}=-sJ_g^{(s)}+\frac12Q_g^{(s)} g,
\qquad
\operatorname{tr}_gJ_g^{(s)}=Q_g^{(s)},
\]
where $J_g^{(s)}$ is the higher $Q$-Ricci tensor.
For $s=2$, it is the $J$-tensor introduced by
Lin and Yuan \cite{LinYuanJTensor2017}.
We obtain the following extension of
Theorem~\ref{Thm:principal}.

\begin{theorem}\label{thm:GJMS-ratio}
Let $2s<n$ and assume that $P_g^s$ is coercive, that its first
eigenvalue is simple, and that a first eigenfunction $u$ is
strictly positive. Then $g$ is Clarke critical for
$\mathcal R_{2s}$ if and only if there exist second
eigenfunctions $v_1,\ldots,v_m$ such that
\[
\psi=\left(\frac{v_1}{u},\ldots,\frac{v_m}{u}\right)
:M\longrightarrow\mathbb S^{m-1}
\]
and, for
\[
g_u=u^{4/(n-2s)}g,
\]
one has
\begin{equation}\label{eq:GJMS-intro-map}
L_{g_u}^{(s),0}\psi
=\bigl(\mathcal R_{2s}(g)-1\bigr)
\frac{n-2s}{2}Q_{g_u}^{(s)}\psi,
\end{equation}
together with the tensorial equation
\begin{equation}\label{eq:GJMS-intro-metric}
H_{g_u}^{(2s)}(\psi)
=\bigl(\mathcal R_{2s}(g)-1\bigr)
\frac{n-2s}{2}E_{g_u}^{Q^{(s)}}.
\end{equation}
\end{theorem}

The right-hand side of~\eqref{eq:GJMS-intro-map}
is normal to the sphere, so that $\psi$ is a critical
point of the reduced GJMS energy under the constraint
$|\psi|=1$.
The tensorial equation is the corresponding higher-order
harmonic--Einstein equation. It can also be written
in the trace-reversed form
\begin{equation}\label{eq:GJMS-intro-trace-reversed}
-\frac1s\left(
H_{g_u}^{(2s)}(\psi)
-\frac{\operatorname{tr}_{g_u}H_{g_u}^{(2s)}(\psi)}{n-2s}g_u
\right)
=\bigl(\mathcal R_{2s}(g)-1\bigr)
\frac{n-2s}{2}J_{g_u}^{(s)}.
\end{equation}
For $s=1$, these equations recover the harmonic--Einstein
system of Theorem~\ref{Thm:principal}.
For $s=2$, they give the corresponding Paneitz system,
with the higher $Q$-Ricci tensor taking the place
of the Ricci tensor.

\subsection{Plan of the paper}

In Section \ref{sec:Clarke}, we recall the necessary material from Clarke analysis
in the framework of Petrides and Tewodrose
\cite{PetridesTewodrose2024}. In particular, we establish the tensorial
characterization of criticality for both the conformal Laplacian
and the GJMS operators.
In Section \ref{sec:L_g} we prove Theorem~\ref{Thm:principal} for the conformal Laplacian, and then we discuss rigidity results
and topological obstructions, together with examples.
In Section \ref{sec:GJMS}, we introduce the curvature tensors and reduced
energies associated with the GJMS operators, and we prove
Theorem~\ref{thm:GJMS-ratio}.

\section{Clarke subdifferential}\label{sec:Clarke}

We start by briefly recalling Clarke criticality, following
Petrides and Tewodrose \cite{PetridesTewodrose2024}.

Let $\mathcal S$ be the space of smooth symmetric two-tensors,
which is the tangent space to $\mathcal M(M)$ at $g$.
If $F$ is a locally Lipschitz functional of the metric and
$h\in\mathcal S$, its generalized upper and lower directional
derivatives at $g$ in the direction $h$ are
\[
F^\circ_+(g;h)
=
\limsup_{\substack{g'\to g\\ t\downarrow0}}
\frac{F(g'+th)-F(g')}{t},
\qquad
F^\circ_-(g;h)
=
\liminf_{\substack{g'\to g\\ t\downarrow0}}
\frac{F(g'+th)-F(g')}{t}.
\]
Here $g'\to g$ is understood in a sufficiently strong
$C^r$-topology, and $t>0$ is small enough so that $g'+th$
is still a Riemannian metric. The Clarke subdifferential
and superdifferential are then defined by
\[
\partial_C^-F(g)
=
\left\{
\ell\in\mathcal S^*:
F^\circ_+(g;h)\geq\ell(h)
\text{ for every }h\in\mathcal S
\right\},
\]
and
\[
\partial_C^+F(g)
=
\left\{
\ell\in\mathcal S^*:
F^\circ_-(g;h)\leq\ell(h)
\text{ for every }h\in\mathcal S
\right\}.
\]
We say that $g$ is Clarke critical for $F$ if
\[
0\in\partial_C^-F(g)\cup\partial_C^+F(g).
\]
See Petrides and Tewodrose
\cite[Definitions~2.11 and~2.14]{PetridesTewodrose2024}
for more details.

If $F$ is strongly differentiable at $g$, both sets reduce
to the ordinary differential $dF_g$. Thus Clarke criticality
is a nonsmooth extension of the usual critical point condition.
This point of view is closely related to the approach developed
by El Soufi and Ilias \cite{ElSoufiIlias2008},
Nadirashvili \cite{Nadirashvili1996Berger}, Petrides and others
for eigenvalue functionals with possible multiplicities.
Along an analytic deformation of the metric, a multiple
eigenvalue splits into finitely many analytic branches,
and criticality is characterized by suitable sign conditions
on the one-sided derivatives of these branches.
The Clarke formalism packages all these possible first
variations into a single convex object, namely the Clarke
sub/superdifferential. This framework is particularly well
suited for quotients and finite combinations of eigenvalues,
since the corresponding Euler--Lagrange equations naturally
take the form of convex balance identities between
first-variation tensors associated with eigenfunctions.

\medskip

We now make this explicit for the ratio
$\lambda_2(M,L_g)/\lambda_1(M,L_g)$. For simplicity, we write
\[
\lambda_i=\lambda_i(M,L_g),
\qquad
E_{\lambda_i}(g)=\ker(L_g-\lambda_i).
\]
If $A$ is a symmetric two-tensor, its trace-free part is
\[
A^0=A-\frac1n(\operatorname{tr}_gA)g.
\]

For a normalized eigenfunction $f\in E_{\lambda_i}(g)$,
the first variation of the corresponding eigenvalue branch
is represented by a linear functional on $\mathcal S$.
We write this functional in tensorial form as
\[
h\longmapsto
\int_M\langle\mathcal T_{\lambda_i}(f),h\rangle_g\,dv_g,
\]
where the explicit expression of the tensor
$\mathcal T_{\lambda_i}(f)$ will be computed
in the next proposition.

Since $\lambda_2>\lambda_1$, Proposition~3.6 of
Petrides and Tewodrose \cite{PetridesTewodrose2024}
identifies $\partial_C^+\lambda_2(M,L_g)$ with the convex
hull of these first-variation functionals.
We use the shorthand
\[
\partial_C\lambda_2(M,L_g)
:=\partial_C^+\lambda_2(M,L_g)
=
\operatorname{conv}
\left\{
\ell_w:
w\in E_{\lambda_2}(g),\
\int_Mw^2\,dv_g=1
\right\}
\subset\mathcal S^*,
\]
where
\[
\ell_w(h)
=
\int_M\langle\mathcal T_{\lambda_2}(w),h\rangle_g\,dv_g.
\]
Thus an element of the Clarke differential of $\lambda_2$
can be written as
\[
\sum_{\alpha=1}^m
c_\alpha\mathcal T_{\lambda_2}(w_\alpha),
\qquad
c_\alpha\geq0,
\qquad
\sum_{\alpha=1}^m c_\alpha=1,
\]
for normalized eigenfunctions
$w_\alpha\in E_{\lambda_2}(g)$.

\subsection{The first variation formula}

We begin by computing the first-variation tensor associated
with the conformal Laplacian under arbitrary metric deformations.

\begin{proposition}
For every normalized eigenfunction $f\in E_\lambda(g)$,
the first-variation tensor introduced above is given explicitly by
\begin{equation}\label{eq:stress-tensor}
\begin{aligned}
\mathcal T_\lambda(f)
={}&-df\otimes df+\frac12|df|_g^2g\\
&+c_n\left(
-f^2\operatorname{Ric}_g
+\nabla^g d(f^2)
+\Delta_g(f^2)g
+\frac12S_gf^2g
\right)
-\frac{\lambda}{2}f^2g.
\end{aligned}
\end{equation}
\end{proposition}

\begin{proof}
Let $\lambda(t)$ be a differentiable eigenvalue branch
along a deformation $g_t$ with $\dot g_0=h$,
with associated eigenfunctions $f_t$ such that $f_0=f$
and $\int_Mf^2\,dv_g=1$.
At $t=0$, the eigenvalue equation gives
\[
\lambda
=
\frac{\displaystyle\int_M
\left(|df|_g^2+c_nS_gf^2\right)\,dv_g}
{\displaystyle\int_Mf^2\,dv_g}.
\]
We need to prove that
\[
\left.\frac d{dt}\right|_{t=0}\lambda(t)
=
\int_M
\langle\mathcal T_\lambda(f),h\rangle_g\,dv_g,
\]
where $\mathcal T_\lambda(f)$ is given
by~\eqref{eq:stress-tensor}.

The terms involving $\dot f_0$ cancel by the eigenvalue
equation. From the variations of the Dirichlet energy
and scalar curvature we have
\[
\left.\frac d{dt}\right|_{t=0}
\int_M|df|_{g_t}^2\,dv_{g_t}
=
\int_M
\left\langle
-df\otimes df+\frac12|df|_g^2g,h
\right\rangle_g\,dv_g,
\]
and
\[
\left.\frac d{dt}\right|_{t=0}S_{g_t}
=
\Delta_g(\operatorname{tr}_gh)
+\operatorname{div}_g\operatorname{div}_gh
-\langle\operatorname{Ric}_g,h\rangle_g.
\]
Therefore, integrating by parts gives
\[
\left.\frac d{dt}\right|_{t=0}
\int_MS_{g_t}f^2\,dv_{g_t}
=
\int_M
\left\langle
-f^2\operatorname{Ric}_g
+\nabla^g d(f^2)
+\Delta_g(f^2)g
+\frac12S_gf^2g,h
\right\rangle_g\,dv_g.
\]
Finally, since
\[
\left.\frac d{dt}\right|_{t=0}
\int_Mf^2\,dv_{g_t}
=
\frac12\int_Mf^2\operatorname{tr}_g(h)\,dv_g,
\]
we deduce that
\[
\left.\frac d{dt}\right|_{t=0}\lambda(t)
=
\int_M
\langle\mathcal T_\lambda(f),h\rangle_g\,dv_g.
\]
\end{proof}

This computation identifies the relevant first-variation
tensors that we will use in the Clarke differential.
The derivatives of the eigenvalue branches issuing from
a multiple eigenvalue are represented by the corresponding
quadratic form on the eigenspace.

For any normalized eigenfunction $f\in E_\lambda(g)$,
since
\[
\operatorname{tr}_g\nabla^g d(f^2)=-\Delta_g(f^2),
\qquad
\Delta_g(f^2)
=2(\lambda-c_nS_g)f^2-2|df|_g^2,
\]
taking the trace of $\mathcal T_\lambda(f)$ we obtain
\begin{align*}
\operatorname{tr}_g\mathcal T_\lambda(f)
&=\frac{n-2}{2}|df|_g^2-\frac n2\lambda f^2
+c_n\left(
\frac{n-2}{2}S_gf^2+(n-1)\Delta_g(f^2)
\right)\\
&=\frac{n-2}{2}|df|_g^2-\frac n2\lambda f^2
+c_n\left(
\frac{n-2}{2}S_gf^2
+2(n-1)(\lambda-c_nS_g)f^2
-2(n-1)|df|_g^2
\right)\\
&=\left(\frac{n-2}{2}-2c_n(n-1)\right)|df|_g^2
+\left(-\frac n2+2c_n(n-1)\right)\lambda f^2\\
&\quad
+c_n\left(\frac{n-2}{2}-2c_n(n-1)\right)S_gf^2\\
&=-\lambda f^2.
\end{align*}
In the same manner, its trace-free part is given by
\begin{equation}\label{eq:TF-stress}
\mathcal T_\lambda(f)^0
=-(df\otimes df)^0
+c_n\left(
\bigl(\nabla^g d(f^2)\bigr)^0
-f^2\operatorname{Ric}_g^0
\right).
\end{equation}

\subsection{Euler--Lagrange equations}
We now derive the Euler--Lagrange characterization of the ratio
functional $\mathcal R_2$. This is essentially contained in Petrides and Tewodrose
\cite[Theorems~3.13 and~3.14]{PetridesTewodrose2024}.
For completeness, we write the full argument.

\begin{proposition}\label{prop: Euler-Lagrange for L_g}
Let $g\in\mathcal M^+(M)$.
Then $g$ is Clarke critical for $\mathcal R_2$
if and only if there exist a first eigenfunction $u>0$
and pairwise $L^2(M,g)$-orthogonal second eigenfunctions
$v_1,\ldots,v_m$ such that
\begin{equation}\label{equa; E-L L_g}
\frac1{\lambda_2}
\sum_{\alpha=1}^m\mathcal T_{\lambda_2}(v_\alpha)
=
\frac1{\lambda_1}\mathcal T_{\lambda_1}(u).
\end{equation}
\end{proposition}

\begin{proof}
Assume first that $g$ is Clarke critical for $\mathcal R_2$.
Since $g\in\mathcal M^+(M)$, the first eigenvalue
$\lambda_1$ is simple.
We choose the positive eigenfunction $u$ normalized by
\[
\int_Mu^2\,dv_g=1.
\]
Thus $\lambda_1$ is $C^1$ in a neighborhood of $g$, and
\[
d\lambda_1(M,L_g)(h)
=
\int_M
\langle\mathcal T_{\lambda_1}(u),h\rangle_g\,dv_g.
\]
By the Clarke description of the differential of
a multiple eigenvalue, an element of
$\partial_C\lambda_2(M,L_g)$ is represented by a convex
combination
\[
\sum_{\alpha=1}^m
c_\alpha\mathcal T_{\lambda_2}(w_\alpha),
\qquad
c_\alpha\geq0,
\qquad
\sum_{\alpha=1}^m c_\alpha=1,
\]
where
\[
w_\alpha\in E_{\lambda_2}(g),
\qquad
\int_Mw_\alpha^2\,dv_g=1.
\]
The chain rule for Clarke differentials applied to
$(x,y)\mapsto y/x$ gives the following necessary
criticality condition:
\[
0=
\frac1{\lambda_1}
\sum_{\alpha=1}^m
c_\alpha\mathcal T_{\lambda_2}(w_\alpha)
-\frac{\lambda_2}{\lambda_1^2}
\mathcal T_{\lambda_1}(u).
\]
Since $f\mapsto\mathcal T_\lambda(f)$ is quadratic,
setting
\[
v_\alpha=\sqrt{c_\alpha}\,w_\alpha
\]
leads to
\[
\mathcal T_{\lambda_2}(v_\alpha)
=
c_\alpha\mathcal T_{\lambda_2}(w_\alpha),
\]
and the desired identity follows.

By Petrides and Tewodrose \cite[Theorem~3.13 (A)]{PetridesTewodrose2024},
the family $(w_\alpha)$ may be chosen
$L^2(M,g)$-orthonormal.
Since $v_\alpha=\sqrt{c_\alpha}\,w_\alpha$, after discarding
the terms with $c_\alpha=0$, the functions
$v_1,\ldots,v_m$ are pairwise $L^2(M,g)$-orthogonal.

Conversely, assume that there exist a first eigenfunction
$u$ and a family of second eigenfunctions $v_1,\ldots,v_m$ satisfying
\[
\frac1{\lambda_2}
\sum_{\alpha=1}^m\mathcal T_{\lambda_2}(v_\alpha)
=
\frac1{\lambda_1}\mathcal T_{\lambda_1}(u).
\]
Since the identity is homogeneous of degree two in the
functions, we may multiply $u,v_1,\ldots,v_m$ by a common positive constant so that $\int_Mu^2\,dv_g=1$.
Taking the trace and using
\[
\operatorname{tr}_g\mathcal T_\lambda(f)
=-\lambda f^2,
\]
we obtain
\[
\sum_{\alpha=1}^m v_\alpha^2=u^2.
\]
Thus, if
\[
b_\alpha=\int_Mv_\alpha^2\,dv_g,
\]
then $b_\alpha\geq0$ and
\[
\sum_{\alpha=1}^m b_\alpha=1.
\]
After discarding the zero terms, set
\[
w_\alpha=\frac{v_\alpha}{\sqrt{b_\alpha}}.
\]
Then
\[
w_\alpha\in E_{\lambda_2}(g),
\qquad
\int_Mw_\alpha^2\,dv_g=1,
\qquad
\mathcal T_{\lambda_2}(v_\alpha)
=
b_\alpha\mathcal T_{\lambda_2}(w_\alpha).
\]
Hence
\[
\eta:=
\sum_{\alpha=1}^m\mathcal T_{\lambda_2}(v_\alpha)
=
\sum_{\alpha=1}^m
b_\alpha\mathcal T_{\lambda_2}(w_\alpha)
\]
is a convex combination of first-variation tensors
associated with normalized eigenfunctions
in $E_{\lambda_2}(g)$.

For $F(x,y)=y/x$, we have
\[
\partial_1F(\lambda_1,\lambda_2)
=-\frac{\lambda_2}{\lambda_1^2},
\qquad
\partial_2F(\lambda_1,\lambda_2)
=\frac1{\lambda_1}>0.
\]
Equation~\eqref{equa; E-L L_g} therefore can be reformulated as
\[
\partial_1F(\lambda_1,\lambda_2)
\mathcal T_{\lambda_1}(u)
+
\partial_2F(\lambda_1,\lambda_2)\eta
=0.
\]
Since $\lambda_1$ is simple and $\lambda_2$ is at the
beginning of its multiplicity block, the positivity of
$\partial_2F(\lambda_1,\lambda_2)$ ensures that the sign
condition in Petrides and Tewodrose
\cite[Theorem~3.14]{PetridesTewodrose2024} is satisfied.
It follows that $g$ is Clarke critical for $\mathcal R_2$.
\end{proof}

\subsection{Euler--Lagrange equations for the GJMS operators}

The same structure extends to every subcritical GJMS operator
$P_g^s$ of order $2s<n$
\cite{GrahamJenneMasonSparling1992}.
The variational theory of their eigenvalues in a fixed
conformal class has recently been developed by Humbert,
Petrides and Premoselli
\cite{HumbertPetridesPremoselli2025},
and fits naturally into the Clarke framework
of Petrides and Tewodrose \cite{PetridesTewodrose2024}.
Since there is no maximum principle in general
for higher-order operators, we assume throughout
this section that
\begin{equation}\label{conditionsGJMS}
\begin{aligned}
\textup{(i)}\quad
& P_g^s \text{ is coercive};\\
\textup{(ii)}\quad
& \text{its first eigenvalue is simple};\\
\textup{(iii)}\quad
& \text{a first eigenfunction can be chosen strictly positive}.
\end{aligned}
\end{equation}

For a fixed smooth function $f$, we define
\[
\mathcal Q_g^{(s)}(f)
=\int_MfP_g^sf\,dv_g,
\qquad
\mathcal N_g(f)=\int_Mf^2\,dv_g.
\]
For an eigenfunction $P_g^sf=\lambda f$,
let $\mathcal T^{(2s)}_\lambda(f)$ denote
the symmetric two-tensor defined, after integration
by parts, by
\begin{equation}\label{eq:GJMS-stress-decomposition}
\int_M
\langle\mathcal T^{(2s)}_\lambda(f),h\rangle_g\,dv_g
=
D_g\mathcal Q_g^{(s)}(f)[h]
-\lambda D_g\mathcal N_g(f)[h],
\end{equation}
for every smooth symmetric two-tensor $h$.
Here $f$ is kept fixed when the metric varies.
If $f$ is $L^2(M,g)$-normalized, the right-hand side
is the first variation of the Rayleigh quotient at $f$.
In particular, if $f$ is the normalized initial eigenfunction
of a differentiable eigenvalue branch $\lambda(t)$
along a deformation $g_t$ with $\dot g_0=h$, then
\begin{equation}\label{eq:GJMS-eigenvalue-stress}
\left.\frac d{dt}\right|_{t=0}\lambda(t)
=
\int_M
\langle\mathcal T^{(2s)}_\lambda(f),h\rangle_g\,dv_g.
\end{equation}
For a multiple eigenvalue, the derivatives of the eigenvalue
branches are the eigenvalues of the corresponding quadratic
form on the eigenspace.
Definition~\eqref{eq:GJMS-stress-decomposition}
also applies to non-normalized eigenfunctions and gives
\[
\mathcal T^{(2s)}_\lambda(cf)
=c^2\mathcal T^{(2s)}_\lambda(f),
\qquad c\in\mathbb R.
\]

We first need its trace.
Let $g_t=e^{2t\omega}g$, where $\omega\in C^\infty(M)$.
From the conformal covariance~\eqref{eq:GJMS-covariance},
\begin{equation}\label{eq:GJMS-quadratic-covariance-conformal}
\int_MfP_{g_t}^sf\,dv_{g_t}
=
\int_M
e^{\frac{n-2s}{2}t\omega}f
P_g^s\!\left(e^{\frac{n-2s}{2}t\omega}f\right)\,dv_g.
\end{equation}
Since $P_g^sf=\lambda f$ and $P_g^s$ is self-adjoint,
\[
\left.\frac d{dt}\right|_{t=0}
\int_MfP_{g_t}^sf\,dv_{g_t}
=
(n-2s)\lambda\int_M\omega f^2\,dv_g.
\]
On the other hand,
\[
\left.\frac d{dt}\right|_{t=0}
\int_Mf^2\,dv_{g_t}
=
n\int_M\omega f^2\,dv_g.
\]
Using~\eqref{eq:GJMS-stress-decomposition}
with $h=2\omega g$, we obtain
\[
2\int_M
\omega\operatorname{tr}_g
\mathcal T^{(2s)}_\lambda(f)\,dv_g
=
-2s\lambda\int_M\omega f^2\,dv_g.
\]
Since $\omega$ is arbitrary, this gives
the pointwise identity
\begin{equation}\label{eq:GJMS-stress-trace}
\operatorname{tr}_g\mathcal T^{(2s)}_\lambda(f)
=-s\lambda f^2.
\end{equation}
For $s=1$, $P_g^1=L_g$ and
$\mathcal T_\lambda^{(2)}(f)=\mathcal T_\lambda(f)$,
so this recovers the trace identity above.
Write
\[
\lambda_i=\lambda_i(M,P_g^s),
\qquad
\mathcal R_{2s}(g)=\frac{\lambda_2}{\lambda_1}.
\]
\begin{proposition}\label{prop:Euler-Lagrange-GJMS}
Let $2s<n$ and assume that conditions \eqref{conditionsGJMS} are satisfied.
Then $g$ is Clarke critical for $\mathcal R_{2s}$
if and only if there exist a first eigenfunction $u>0$
and second eigenfunctions $v_1,\ldots,v_m$ such that
\begin{equation}\label{eq:GJMS-Clarke-balance}
\frac1{\lambda_2}\sum_{\alpha=1}^m
\mathcal T^{(2s)}_{\lambda_2}(v_\alpha)
=
\frac1{\lambda_1}
\mathcal T^{(2s)}_{\lambda_1}(u).
\end{equation}
\end{proposition}

\begin{proof}
We apply the same characterization of Petrides
and Tewodrose
\cite[Theorems~3.13 and~3.14]{PetridesTewodrose2024}
as in Proposition~\ref{prop: Euler-Lagrange for L_g},
with numerator $\mathcal Q_g^{(s)}$ and denominator
$\mathcal N_g$ on $H^s(M)$.
Their abstract assumptions hold locally in a sufficiently
strong $C^r$-topology on the metrics: the quadratic forms
depend smoothly on the metric, coercivity persists near $g$,
and $H^s(M)\hookrightarrow L^2(M)$ is compact.

Choose $u>0$ with $\int_Mu^2\,dv_g=1$.
For $F(x,y)=y/x$, the partial derivatives
at $(\lambda_1,\lambda_2)$ are
\[
\partial_1F=-\frac{\lambda_2}{\lambda_1^2},
\qquad
\partial_2F=\frac1{\lambda_1}.
\]
The necessary condition therefore gives normalized
second eigenfunctions $w_\alpha$ and coefficients
$c_\alpha\geq0$, with $\sum_\alpha c_\alpha=1$,
such that
\[
0=
\frac1{\lambda_1}\sum_\alpha
c_\alpha\mathcal T^{(2s)}_{\lambda_2}(w_\alpha)
-\frac{\lambda_2}{\lambda_1^2}
\mathcal T^{(2s)}_{\lambda_1}(u).
\]
Setting $v_\alpha=\sqrt{c_\alpha}\,w_\alpha$
gives~\eqref{eq:GJMS-Clarke-balance}.

Conversely, suppose that~\eqref{eq:GJMS-Clarke-balance}
holds. Since the identity is homogeneous of degree two in the
functions, we may multiply $u,v_1,\ldots,v_m$ by a common
positive constant so that
\[
\int_M u^2\,dv_g=1.
\]
Taking the trace and using~\eqref{eq:GJMS-stress-trace}, we obtain
\[
-s\sum_{\alpha=1}^m v_\alpha^2=-su^2.
\]
Integrating over $M$ therefore gives
\[
b_\alpha:=\int_Mv_\alpha^2\,dv_g,
\qquad
\sum_{\alpha=1}^m b_\alpha=1.
\]
After discarding zero terms, the functions
$w_\alpha=v_\alpha/\sqrt{b_\alpha}$
are normalized second eigenfunctions, and the tensor balance
becomes
\[
0=
-\frac{\lambda_2}{\lambda_1^2}
\mathcal T^{(2s)}_{\lambda_1}(u)
+\frac1{\lambda_1}\sum_{\alpha=1}^m
b_\alpha\mathcal T^{(2s)}_{\lambda_2}(w_\alpha).
\]
Since $\lambda_1$ is simple and $\lambda_2>\lambda_1$,
the sufficient condition of Petrides and Tewodrose
\cite[Theorem~3.14]{PetridesTewodrose2024}
applies to the same function $F(x,y)=y/x$:
the total weights on the first and second eigenspaces
are respectively $-\lambda_2/\lambda_1^2$
and $1/\lambda_1$.
It yields $0\in\partial_C^+\mathcal R_{2s}(g)$.
Thus $g$ is Clarke critical for $\mathcal R_{2s}$.
\end{proof}

\section{Proof of the Main result and Consequences }\label{sec:L_g}
\subsection{Proof of Theorem \ref{Thm:principal}}
 The proof of Theorem \ref{Thm:principal} is based on a simple but useful decomposition of the Euler--Lagrange
equation into its trace and trace-free parts. The trace gives the spherical constraint and the harmonicity of the eigenvalue ratio map, while the trace-free part, after the conformal change \(g_u=u^{4/(n-2)}g\), identifies the pull-back metric of this map with the Ricci tensor of \(g_u\).

\begin{proof}[Proof of Theorem \ref{Thm:principal}]
Assume that $g$ is Clarke critical for $\mathcal R_2$.
By Proposition~\ref{prop: Euler-Lagrange for L_g}, choose
a positive first eigenfunction $u$, normalized by
$\int_M u^2\,dv_g=1$, and pairwise
$L^2(M,g)$-orthogonal second eigenfunctions
$v_1,\ldots,v_m$ such that
$$
\frac1{\lambda_2}\sum_{\alpha=1}^m\mathcal T_{\lambda_2}(v_\alpha)
=
\frac1{\lambda_1}\mathcal T_{\lambda_1}(u).
$$ 

Using the identity \(\operatorname{tr}_g\mathcal T_\lambda(f)=-\lambda f^2\), we obtain
$$
\sum_{\alpha=1}^m v_\alpha^2=u^2.
$$
Hence $\psi:=(v_1/u,\ldots,v_m/u)$ maps $M$ into $\mathbb S^{m-1}$.

From the conformal covariance of the conformal Laplacian, since
$$
L_{g_u}\varphi
=u^{-\frac{n+2}{n-2}}L_g(u\varphi), \quad g_u=u^{4/(n-2)}g,
$$
we have
$$
c_nS_{g_u}=\lambda_1u^{-\frac4{n-2}},
\qquad
L_{g_u}\psi=\lambda_2u^{-\frac4{n-2}}\psi.
$$
Therefore
$$
\Delta_{g_u}\psi
=(\lambda_2-\lambda_1)u^{-\frac4{n-2}}\psi.
$$
Since $\sum_\alpha\psi_\alpha^2=1$, and $\frac12\Delta_{g_u}(\psi_\alpha^2)
=\psi_\alpha\Delta_{g_u}\psi_\alpha-|d\psi_\alpha|_{g_u}^2,$
we obtain
$$
|d\psi|_{g_u}^2=(\lambda_2-\lambda_1)u^{-\frac4{n-2}},
$$
hence
$$
\Delta_{g_u}\psi=|d\psi|_{g_u}^2\psi.
$$
This proves that the map $\psi$ is harmonic into the sphere $\mathbb{S}^{m-1}$ in the metric $g_u$.

It remains to prove the tensor identity. We have
$$
\psi^*g_{\mathbb S^{m-1}}=\sum_{\alpha=1}^m d\psi_\alpha\otimes d\psi_\alpha,
$$
and by $v_\alpha=u\psi_\alpha$ and $\sum_\alpha\psi_\alpha^2=1$, it is straightforward to see that
\begin{equation}\label{eq:dv-sum}
\sum_{\alpha=1}^m dv_\alpha\otimes dv_\alpha
=du\otimes du+u^2\psi^*g_{\mathbb S^{m-1}},
\qquad
\sum_{\alpha=1}^m\nabla^g d(v_\alpha^2)
=\nabla^g d(u^2).
\end{equation}
Taking trace-free parts in \eqref{equa; E-L L_g}, using \eqref{eq:TF-stress} and \eqref{eq:dv-sum},
\begin{equation}\label{eq:TF-before-conformal}
\lambda_1u^2(\psi^*g_{\mathbb S^{m-1}})^0
=(\lambda_2-\lambda_1)
\left[
(du\otimes du)^0
-c_n\left(
\bigl(\nabla^g d(u^2)\bigr)^0
-u^2\operatorname{Ric}_g^0
\right)
\right].
\end{equation}
The expression on the right-hand side of \eqref{eq:TF-before-conformal} is none other than the trace-free Ricci tensor of $g_u=u^{4/(n-2)}g$, 
\begin{equation}\label{eq:conformal-ricci-identity}
u^{-2}\left[
(du\otimes du)^0
-c_n\left(
\bigl(\nabla^g d(u^2)\bigr)^0
-u^2\operatorname{Ric}_g^0
\right)
\right]
=c_n\operatorname{Ric}_{g_u}^0.
\end{equation}
Thus, we deduce that
$$
(\psi^*g_{\mathbb S^{m-1}})^0=c_n\left(\frac{\lambda_2}{\lambda_1}-1\right)\operatorname{Ric}_{g_u}^0.
$$
To complete the proof of this part, notice that the traces with respect to $g_u$ also coincide,
$$
\operatorname{tr}_{g_u}\psi^*g_{\mathbb S^{m-1}}=|d\psi|_{g_u}^2=(\lambda_2-\lambda_1)u^{-\frac4{n-2}}, \quad \text{and} \quad c_nS_{g_u}=\lambda_1u^{-\frac4{n-2}}.
$$
Hence
$$
\operatorname{tr}_{g_u}\psi^*g_{\mathbb S^{m-1}}
=c_n\left(\frac{\lambda_2}{\lambda_1}-1\right)S_{g_u}.
$$
This finally implies that
$$
\psi^*g_{\mathbb S^{m-1}}
=
 c_n\left(\frac{\lambda_2}{\lambda_1}-1\right)\operatorname{Ric}_{g_u}.
$$

Conversely, assume that there exist eigenfunctions \(u>0\),
\(v_1,\ldots,v_m\) satisfying
\[
L_g u=\lambda_1 u,\qquad L_g v_\alpha=\lambda_2 v_\alpha,
\]
such that the map
\[
\psi=\left(\frac{v_1}{u},\ldots,\frac{v_m}{u}\right):M\longrightarrow
\mathbb S^{m-1}
\]
is harmonic with respect to \(g_u=u^{4/(n-2)}g\) and satisfies
\[
\psi^*g_{\mathbb S^{m-1}}
=
c_n\left(\frac{\lambda_2}{\lambda_1}-1\right)\operatorname{Ric}_{g_u}.
\]
Since \(|\psi|^2=1\), we have $\sum_{\alpha=1}^m v_\alpha^2=u^2$.
From
$\operatorname{tr}_g \mathcal T_\lambda(f)=-\lambda f^2,
$ we obtain
\[
\operatorname{tr}_g\left(
\frac1{\lambda_2}\sum_{\alpha=1}^m \mathcal T_{\lambda_2}(v_\alpha)
\right)
=
-u^2
=
\operatorname{tr}_g\left(
\frac1{\lambda_1}\mathcal T_{\lambda_1}(u)
\right).
\]
Thus the traces agree.

Taking the trace-free part of the harmonic--Einstein equation, we recover exactly
\[
\lambda_1u^2
\left(\psi^*g_{\mathbb S^{m-1}}\right)^0
=
(\lambda_2-\lambda_1)
\left[
(du\otimes du)^0
-
c_n\left(\bigl(\nabla^g d(u^2)\bigr)^0-u^2\operatorname{Ric}_g^0\right)
\right].
\]
This last identity is precisely  the trace-free part of
\[
\frac1{\lambda_2}\sum_{\alpha=1}^m \mathcal T_{\lambda_2}(v_\alpha)
=
\frac1{\lambda_1}\mathcal T_{\lambda_1}(u).
\]
Since the traces also agree, we obtain the equality of the full tensor. We conclude by Proposition \ref{prop: Euler-Lagrange for L_g} to deduce that \(g\) is Clarke critical for \(\mathcal R_2\).

The proof of Theorem \ref{Thm:principal} is complete.
\end{proof}

\subsection{Rigidity, topology and examples}

As consequences, we obtain the following rigidity statements.

\begin{corollary}\label{cor:spherical-rigidity-intro}
If $g$ is critical for $\mathcal R_2$, then
$$
\frac{\lambda_2(M,L_g)}{\lambda_1(M,L_g)}\geq \frac{n+2}{n-2}.
$$
Moreover, if equality holds, then $(M,g_u)$ is a round sphere.
\end{corollary}
Notice that the ratio $\frac{n+2}{n-2}$ is exactly
the value of $\frac{\lambda_2(L_{g_0})}{\lambda_1(L_{g_0})}$
for the round sphere $(\mathbb S^n,g_0)$.
In fact, the identity map on $\mathbb S^n$ is harmonic, the constant function $1$ and the coordinate functions
are respectively first and second eigenfunctions of $L_{g_0}$.
Taking $u\equiv1$ and $\psi=\operatorname{id}_{\mathbb S^n}$,
we easily see that equation~\eqref{eq:Einstein-harmonic}
is satisfied, and therefore $g_0$ is critical for $\mathcal R_2$
by Theorem~\ref{Thm:principal}.

\begin{proof}[Proof of Corollary \ref{cor:spherical-rigidity-intro}]
Assume that $g$ is critical, and let $g_u=u^\frac{4}{n-2}g$, where $u>0$ is the first positive eigenfunction of $g$. First we take the $g_u$-trace in \eqref{eq:Einstein-harmonic} to obtain
\begin{equation}\label{eq:energy-scalar}
|d\psi|_{g_u}^2=aS_{g_u}, \quad \text{where} \quad a:=c_n\left(\frac{\lambda_2}{\lambda_1}-1\right).
\end{equation}
We use the Bochner identity for sphere-valued  harmonic maps,
$$
\frac12\Delta_{g_u}|d\psi|_{g_u}^2
=-|\nabla^{g_u}d\psi|_{g_u}^2
-\langle \operatorname{Ric}_{g_u},\psi^*g_{\mathbb S^{m-1}}\rangle_{g_u}
+|d\psi|_{g_u}^4
-|\psi^*g_{\mathbb S^{m-1}}|_{g_u}^2,
$$
where $\nabla^{g_u}d\psi$ denotes the intrinsic
covariant derivative of $d\psi$, computed using
the connections of $g_u$ and $\mathbb S^{m-1}$.
Integrating this identity and using
$\psi^*g_{\mathbb S^{m-1}}
=a\operatorname{Ric}_{g_u}$ 
together with \eqref{eq:energy-scalar}, we obtain
\begin{equation}\label{eq:Bochner-identity}
\int_M |\nabla^{g_u}d\psi|_{g_u}^2\,dv_{g_u}
+a(1+a)\int_M |\operatorname{Ric}_{g_u}|_{g_u}^2\,dv_{g_u}
=
 a^2\int_M S_{g_u}^2\,dv_{g_u}.
\end{equation}
On the other hand, since 
$$
|\operatorname{Ric}_{g_u}|_{g_u}^2
=|\operatorname{Ric}_{g_u}^0|_{g_u}^2+\frac1nS_{g_u}^2,
$$
where $\operatorname{Ric}_{g_u}^0$ is the trace-free part of $\operatorname{Ric}_{g_u}$,
we deduce that
\begin{equation}\label{eq:Bochner-gap}
\int_M |\nabla^{g_u}d\psi|_{g_u}^2\,dv_{g_u}
+a(1+a)\int_M |\operatorname{Ric}_{g_u}^0|_{g_u}^2\,dv_{g_u}
=
 a\left(a-\frac{1+a}{n}\right)
\int_M S_{g_u}^2\,dv_{g_u},
\end{equation}
The left-hand side is non-negative. Moreover, by the
conformal covariance of $L_g$,
\[
c_n S_{g_u}= 
L_{g_u}(1)
=u^{-\frac{n+2}{n-2}}L_g(u)
=\lambda_1(M,L_g)u^{-4/(n-2)}>0.
\]
Thus $S_{g_u}>0$, and we deduce that
$$
a\ge\frac1{n-1},
$$
which is exactly the claimed inequality
$$
\frac{\lambda_2(M,L_g)}{\lambda_1(M,L_g)}\ge\frac{n+2}{n-2}.
$$

If equality holds, then $a=1/(n-1)$ and the right-hand side of \eqref{eq:Bochner-gap} vanishes. Hence $\nabla^{g_u}d\psi=0$ and $\operatorname{Ric}_{g_u}^0=0$. Thus $g_u$ is Einstein with positive scalar curvature,
$$
\operatorname{Ric}_{g_u}=(n-1)Kg_u
$$
for some $K>0$. From \eqref{eq:energy-scalar}, we deduce that $|d\psi|_{g_u}^2=nK$, and therefore there exists a non-constant coordinate function $\psi_\alpha$ that satisfies
$$
\Delta_{g_u}\psi_\alpha=nK\psi_\alpha.
$$
By Lichnerowicz's theorem, since
$\operatorname{Ric}_{g_u}=(n-1)Kg_u$, we obtain
$\lambda_1(M,\Delta_{g_u}) \geq nK$.
Together with the fact that $\psi_\alpha$ is a nonconstant eigenfunction
with eigenvalue $nK$, we deduce that $\lambda_1(M,\Delta_{g_u}) = nK$.
Obata's rigidity theorem therefore implies that
$(M,g_u)$ is isometric to the round sphere
of sectional curvature $K$.
\end{proof}

The following refinement is motivated by the rank information contained in the harmonic--Einstein equation. Indeed, since
\[
\psi^*g_{\mathbb S^{m-1}}=a\,\operatorname{Ric}_{g_u},
\]
the rank of \(\operatorname{Ric}_{g_u}\) is exactly the rank of \(d\psi\). Thus the number of positive Ricci directions measures the effective dimension of the spherical part detected by the map $\psi$. In particular, when \(\operatorname{Ric}_{g_u}\) is positive definite, the map \(\psi\) is an immersion and one necessarily has \(\dim E_{\lambda_2}(g)\ge n+1\). The next corollary shows that, more generally, the maximal rank \(r\) of \(\operatorname{Ric}_{g_u}\) controls the eigenvalue ratio and that equality forces a spherical splitting in precisely those \(r\) directions.

\begin{corollary}\label{cor:rank-r-intro}
Let $g\in\mathcal M^+(M)$ be critical for $\mathcal R_2$, and let
$$
r=\max_{x\in M}\operatorname{rank}\operatorname{Ric}_{g_u}(x),
\qquad
 g_u=u^{4/(n-2)}g,
$$
where $u>0$ is a first eigenfunction of $L_g$.
Then $r\ge 2$, and
$$
\frac{\lambda_2(M,L_g)}{\lambda_1(M,L_g)}
\ge
1+\frac{4(n-1)}{(n-2)(r-1)}.
$$
Moreover, if equality holds, then the universal cover of $(M,g_u)$ is isometric to
$$
\mathbb S_K^r\times (F^{n-r},h_F),
$$
where $\mathbb S_K^r$ is the round $r$-sphere of sectional curvature $K>0$, and $(F^{n-r},h_F)$ is complete, simply connected, and Ricci-flat.
\end{corollary}

\begin{proof}[Proof of Corollary \ref{cor:rank-r-intro}]
    
Assume that $g$ is critical. Let
$$
r=\max_{x\in M}\operatorname{rank}\operatorname{Ric}_{g_u}(x).
$$
Since $\operatorname{Ric}_{g_u}\ge0$, at every point we have
\begin{equation}\label{eq:rank-cauchy}
|\operatorname{Ric}_{g_u}|_{g_u}^2\ge\frac{S_{g_u}^2}{r}.
\end{equation}
Indeed, if $\rho_1,\ldots,\rho_k$ are the non-zero eigenvalues of $\operatorname{Ric}_{g_u}$ and $k\le r$, then
$$
S_{g_u}^2=\left(\sum_{i=1}^k\rho_i\right)^2
\le r\sum_{i=1}^k\rho_i^2
=r|\operatorname{Ric}_{g_u}|_{g_u}^2.
$$
Substituting~\eqref{eq:rank-cauchy}
into~\eqref{eq:Bochner-identity} and using
the nonnegativity of $|\nabla^{g_u}d\psi|_{g_u}^2$, we obtain
\[
a^2\int_M S_{g_u}^2\,dv_{g_u}
\geq
\frac{a(1+a)}{r}
\int_M S_{g_u}^2\,dv_{g_u}.
\]
Since $S_{g_u}>0$ and $a>0$, we obtain $a(r-1)\geq1$. Thus $r\geq2$ and $a\geq1/(r-1)$, which gives 
$$
\frac{\lambda_2(M,L_g)}{\lambda_1(M,L_g)}
\ge
1+\frac1{c_n(r-1)}
=
1+\frac{4(n-1)}{(n-2)(r-1)}.
$$

Assume that equality holds. Then all inequalities above are equalities. Hence $\nabla^{g_u}d\psi=0,$ and equality holds in
$|\operatorname{Ric}_{g_u}|_{g_u}^2\ge \frac{S_{g_u}^2}{r}.$

Therefore $\operatorname{Ric}_{g_u}$ has rank exactly $r$ everywhere and all its non-zero eigenvalues coincide. At each point $x\in M$, the nonzero eigenvalues are therefore
all equal to $\mu(x)=S_{g_u}(x)/r>0$.

Since $\psi^*g_{\mathbb S^{m-1}}=a\,\operatorname{Ric}_{g_u},$
and since \(\nabla^{g_u}d\psi=0\), the tensor \(\psi^*g_{\mathbb S^{m-1}}\) is parallel, which implies that
$$
\nabla^{g_u}\operatorname{Ric}_{g_u}=0.
$$
It follows that the positive eigenspace of \(\operatorname{Ric}_{g_u}\) and its kernel are parallel distributions. By the de Rham decomposition theorem, the universal Riemannian cover of \((M,g_u)\) splits isometrically as
$
(\widetilde M,\widetilde g_u)\cong (N^r,g_N)\times(F^{n-r},h_F)
$
where
$$
\operatorname{Ric}_N=\mu g_N,
\qquad
\operatorname{Ric}_F=0.
$$

Write \(\mu=(r-1)K\) with \(K>0\). Since \(\operatorname{Ric}_{g_u}=0\) along the factor \(F\), we have
$$
\psi^*g_{\mathbb S^{m-1}}=0
$$
on vectors tangent to \(F\). Thus \(\psi\) is constant along \(F\), and therefore only depends on the \(N\)-factor.
Moreover, since equality holds, we have \(a=1/(r-1)\). Hence
$$
\psi^*g_{\mathbb S^{m-1}}
=
a\,\operatorname{Ric}_N
=
K g_N.
$$
Taking the trace on \(N\), we obtain
$$
|d\psi|_{g_N}^2=rK.
$$
Since \(\psi\) is harmonic into the sphere,
$$
\Delta_N\psi=|d\psi|_{g_N}^2\psi,
$$
there exists some non-constant coordinate function \(\psi_\alpha\) that satisfies
$$
\Delta_N\psi_\alpha=rK\psi_\alpha.
$$

On the other hand, we also have
$$
\operatorname{Ric}_N=(r-1)Kg_N,
$$
so by the same argument as in the proof of Corollary \ref{cor:spherical-rigidity-intro}, the Lichnerowicz theorem implies that \(N\) is isometric to the round sphere \(\mathbb S_K^r\). Hence the universal cover of \((M,g_u)\) is isometric to
$$
\mathbb S_K^r\times(F^{n-r},h_F),
$$
where $F$ is Ricci-flat.
\end{proof}

\begin{remark}
Let \(m_2=\dim E_{\lambda_2}(g)\). Since the map \(\psi\) is defined by ratios of second eigenfunctions, it takes values in \(\mathbb S^{m_2-1}\). Moreover, from
\[
\psi^*g_{\mathbb S^{m_2-1}}=a\,\operatorname{Ric}_{g_u},
\]
so that
\[
\operatorname{rank}d\psi
=
\operatorname{rank}\operatorname{Ric}_{g_u}.
\]
At a point where the rank r of $Ric_{g_u}$ is maximal, we therefore obtain
\[
r\le m_2-1.
\]
Since Corollary \ref{cor:rank-r-intro} gives \(r\ge2\), it follows that
\[
\dim E_{\lambda_2}(g)\ge r+1\ge3.
\]
\end{remark}

The rank estimate first gives a direct obstruction
to the existence of Yamabe-positive metrics critical
for $\mathcal R_2$, expressed in terms of the first
Betti number. Indeed, since
$
\operatorname{Ric}_{g_u}\geq0,$
the Bochner formula for harmonic one-forms shows that every harmonic
one-form on \((M,g_u)\) is parallel and that its metric dual belongs to
the kernel of \(\operatorname{Ric}_{g_u}\).

More precisely, if \(\omega\) is harmonic, then by the Bochner--Weitzenb\"ock formula we have
\[
0
=
\int_M
\left(
|\nabla\omega|_{g_u}^2
+
\operatorname{Ric}_{g_u}
(\omega^\sharp,\omega^\sharp)
\right)
\,dv_{g_u}.
\]
Both terms are non-negative, and hence
\[
\nabla\omega=0,
\qquad
\omega^\sharp\in\ker\operatorname{Ric}_{g_u}.
\]
A basis of harmonic one-forms is pointwise linearly independent because
all these forms are parallel. Thus,
\[
b_1(M)
\leq
\dim\ker\operatorname{Ric}_{g_u}(x)
\]
at every point \(x\in M\). We deduce that
\[
r\leq n-b_1(M).
\]
Since Corollary~\ref{cor:rank-r-intro} gives \(r\geq2\), we obtain the bound
\[
b_1(M)\leq n-2.
\]

This argument only uses harmonic one-forms on \(M\) itself. Applying the
same argument to finite coverings and using the Cheeger--Gromoll
structure theorem gives a substantially stronger obstruction involving
the whole fundamental group.

Recall that a group \(G\) is called \emph{virtually abelian} if it
contains an abelian subgroup of finite index. If \(G\) is finitely
generated and virtually abelian, then it contains a finite-index
subgroup isomorphic to \(\mathbb Z^q\) for some \(q\geq0\). The integer
\(q\) is called the \emph{virtual abelian rank} of \(G\).

\begin{corollary}
\label{cor:virtual-rank-intro}
Let \(g\in\mathcal M^+(M)\) be critical for \(\mathcal R_2\), and let
\[
g_u=u^{4/(n-2)}g,
\qquad
r=\max_{x\in M}\operatorname{rank}\operatorname{Ric}_{g_u}(x).
\]
Then \(\pi_1(M)\) is virtually abelian, and its virtual abelian rank \(q\) satisfies
\begin{equation}\label{eq:Betti-virtual-rank-chain}
b_1(M)\leq q\leq n-r\leq n-2
\end{equation}
and
\begin{equation}\label{eq:virtual-rank-bound}
\frac{\lambda_2(M,L_g)}{\lambda_1(M,L_g)}
\geq
1+\frac{4(n-1)}{(n-2)(n-q-1)}.
\end{equation}
Moreover, if equality holds in \eqref{eq:virtual-rank-bound}, then \(u\) is
constant and
\[
(\widetilde M,\widetilde g_u)\cong S_K^{\,n-q}\times\mathbb R^q
\]
for some \(K>0\).

The bound is sharp for every \(q\in\{0,\ldots,n-2\}\). More
precisely, set \(r=n-q\), let \(\mathbb S_K^r\) be the round \(r\)-sphere of
sectional curvature \(K>0\), and let \((T^q,h)\) be a flat torus such
that
\[
\lambda_1(T^q,\Delta_h)\geq rK.
\]
Then the product metric \(g=g_{\mathbb S_K^r}+h\) on \(\mathbb S^r\times T^q\) is critical
for \(\mathcal R_2\), and equality holds in \eqref{eq:virtual-rank-bound}.
\end{corollary}

\begin{proof}
Assume that \(g\) is critical for \(\mathcal R_2\), and let \(u\),
\(v_1,\ldots,v_m\) and \(\psi\) be given by Theorem~\ref{Thm:principal}.
Since \(\lambda_2(M,L_g)>\lambda_1(M,L_g)>0\), the coefficient in
\eqref{eq:Einstein-harmonic} is positive, hence
\[
\operatorname{Ric}_{g_u}\geq0.
\]
By the Cheeger--Gromoll structure theorem \cite{CheegerGromoll1971},
\(\pi_1(M)\) is virtually abelian. Let \(q\) be its virtual abelian rank
and choose a connected finite covering \(p:\widehat M\to M\) such that
\(\pi_1(\widehat M)\simeq\mathbb Z^q\). Set \(\widehat g_u=p^*g_u\). Then
\(\widehat M\) is closed, \(\operatorname{Ric}_{\widehat g_u}\geq0\) and
\(b_1(\widehat M)=q\). Since \(p\) is a local isometry, pulling back
harmonic one-forms is injective, so
\[
b_1(M)\leq b_1(\widehat M)=q.
\]

We now compare \(q\) with \(r\). Let \(\omega\) be a harmonic one-form on
\((\widehat M,\widehat g_u)\). Integrating the Bochner--Weitzenb\"ock formula, we obtain
\[
0=\int_{\widehat M}
\Bigl(|\nabla\omega|^2_{\widehat g_u}
+\operatorname{Ric}_{\widehat g_u}(\omega^\sharp,\omega^\sharp)\Bigr)
dv_{\widehat g_u},
\]
and both terms are nonnegative. Hence \(\omega\) is parallel and
\(\omega^\sharp\in\ker\operatorname{Ric}_{\widehat g_u}\). Let
\(\omega_1,\ldots,\omega_q\) be a basis of harmonic one-forms. They are
pointwise linearly independent: a linear combination which vanishes at
one point is parallel, hence vanishes identically, and so all its
coefficients are zero. Therefore, at every \(x\in\widehat M\),
\[
q\leq\dim\ker\operatorname{Ric}_{\widehat g_u}(x)
=n-\operatorname{rank}\operatorname{Ric}_{\widehat g_u}(x).
\]
Since \(p\) is a local isometry, this gives \(q\leq n-r\), and
Corollary~\ref{cor:rank-r-intro} yields \(r\geq2\). This proves the first point
\eqref{eq:Betti-virtual-rank-chain}.

By Corollary~\ref{cor:rank-r-intro} and since \(r\leq n-q\),
\[
\frac{\lambda_2(M,L_g)}{\lambda_1(M,L_g)}
\geq1+\frac{4(n-1)}{(n-2)(r-1)}
\geq1+\frac{4(n-1)}{(n-2)(n-q-1)},
\]
which is \eqref{eq:virtual-rank-bound}.

Assume now that equality holds in \eqref{eq:virtual-rank-bound}. Then both
inequalities above are equalities. In particular \(r=n-q\), and equality
holds in Corollary~\ref{cor:rank-r-intro}, so
\[
(\widetilde M,\widetilde g_u)\simeq \mathbb S_K^r\times(F^q,h_F),
\]
where \(K>0\) and \((F^q,h_F)\) is complete, simply connected and
Ricci-flat. The scalar curvature \(S_{g_u}=r(r-1)K\) is then constant.
By conformal covariance,
\[
c_nS_{g_u}=L_{g_u}(1)=u^{-\frac{n+2}{n-2}}L_g(u)
=\lambda_1(M,L_g)\,u^{-\frac4{n-2}},
\]
and since \(S_{g_u}\) and \(\lambda_1(M,L_g)\) are positive constants,
\(u\) is constant.

It remains to identify \(F\). The parallel fields
\(\omega_1^\sharp,\ldots,\omega_q^\sharp\) lift to pointwise independent
parallel fields on \(\widetilde M\), which lie in
\(\ker\operatorname{Ric}_{\widetilde g_u}\). This kernel is exactly \(TF\),
because the spherical factor has positive Ricci curvature. Thus \(TF\)
carries a parallel frame, so \(F\) is flat, and being complete and simply
connected it is isometric to \(\mathbb R^q\) (a point if \(q=0\)).

Finally, we prove sharpness. Let \(q\in\{0,\ldots,n-2\}\), \(r=n-q\geq2\)
and \(K>0\). For \(q\geq1\), let \((T^q,h)\) be a flat torus with
\(\lambda_1(T^q,\Delta_h)\geq rK\); this is obtained by scaling any flat
metric by a small positive constant. For \(q=0\), \(T^0\) is a point and
there is no condition. On \(g=g_{\mathbb S_K^r}+h\) we have \(S_g=r(r-1)K>0\) and
\(L_g=\Delta_g+c_nr(r-1)K\), so
\[
\lambda_1(M,L_g)=c_nr(r-1)K,\qquad u\equiv1.
\]
The spectrum of \(\Delta_g\) consists of sums of eigenvalues of the
factors. Since the first nonzero eigenvalue of \(\mathbb S_K^r\) is \(rK\), the
assumption on the torus leads to \(\lambda_1(M,\Delta_g)=rK\), hence
\(\lambda_2(M,L_g)=c_nr(r-1)K+rK\) and
\[
\frac{\lambda_2(M,L_g)}{\lambda_1(M,L_g)}
=1+\frac1{c_n(r-1)}
=1+\frac{4(n-1)}{(n-2)(n-q-1)}.
\]
To check criticality, realize \(\mathbb S_K^r\) as the sphere of radius
\(K^{-1/2}\) in \(\mathbb R^{r+1}\) and set
\(\psi(x,y)=\sqrt K\,x\in\mathbb S^r\). Its components satisfy
\(\Delta_g\psi_\alpha=rK\psi_\alpha\), so they are second eigenfunctions
of \(L_g\). Moreover \(|d\psi|_g^2=rK\) and
\(\Delta_g\psi=|d\psi|_g^2\psi\), so \(\psi\) is harmonic. Since
\[
\psi^*g_{\mathbb S^r}=Kg_{\mathbb S_K^r},
\qquad
\operatorname{Ric}_g=(r-1)Kg_{\mathbb S_K^r}\oplus0,
\qquad
c_n\Bigl(\frac{\lambda_2}{\lambda_1}-1\Bigr)=\frac1{r-1},
\]
equation \eqref{eq:Einstein-harmonic} of Theorem~\ref{Thm:principal} holds with \(g_u=g\). Hence \(g\) is critical for \(\mathcal R_2\).
Since \(r\geq2\), \(\pi_1(\mathbb S^r\times T^q)\simeq\mathbb Z^q\), so the virtual
abelian rank is exactly \(q\).
\end{proof}

The inequality \eqref{eq:virtual-rank-bound} shows that the fundamental
group quantitatively constrains the eigenvalue ratio of a critical metric.
For \(q=0\) it reduces to the spherical estimate
\(\frac{\lambda_2}{\lambda_1}\geq\frac{n+2}{n-2}\), and the bound becomes
strictly stronger as \(q\) increases. The products \(S_K^{n-q}\times T^q\)
show that it is optimal for every rank allowed by
\eqref{eq:Betti-virtual-rank-chain}.

The obstruction coming from the whole fundamental group is also strictly
stronger than the one coming from \(b_1(M)\) alone. This is illustrated by
the following dichotomy.

\begin{corollary}
\label{cor:connected-sum-dichotomy-intro}
Let \(n\geq3\) and \(k\geq1\). Then
\[
\#_k(\mathbb S^{n-1}\times\mathbb S^1)
\]
admits a critical metric for \(\mathcal R_2\) if and only if \(k=1\).

For every \(k\geq2\), this manifold admits metrics in
\(\mathcal M^+(M)\), but none of them is critical for \(\mathcal R_2\).
\end{corollary}

Indeed, we have
\[
\pi_1\left(
\#_k(\mathbb S^{n-1}\times\mathbb S^1)
\right)
\cong F_k,
\]
where \(F_k\) is the free group on \(k\) generators. For \(k\geq2\),
the group \(F_k\) is not virtually abelian, whereas the case \(k=1\) is
realized by the critical product metrics on \(S^{n-1}\times S^1\).

The product of the round metric on $\mathbb S^{n-1}$ with any metric on
$\mathbb S^1$ has positive scalar curvature.
By Gromov and Lawson \cite[Theorem~A]{GromovLawson1980}, this property
is preserved under connected sums in dimension $n\geq3$.
Hence $\#_k(\mathbb S^{n-1}\times\mathbb S^1)$ admits a metric of positive
scalar curvature for every $k\geq1$, and therefore its space
$\mathcal M^+$ is nonempty.

\begin{corollary}
\label{cor:further-topological-obstructions-intro}
Let \(M^n\) be a closed manifold admitting a critical metric for
\(\mathcal R_2\). Then:
\begin{enumerate}[label=\textup{(\roman*)}]
\item \(M\) is not aspherical;

\item if \(\chi(M)\neq0\), then \(\pi_1(M)\) is finite;

\item if \(n\geq4\) and \(\Sigma_\gamma\) is a closed orientable surface
of genus \(\gamma\geq2\), then
\[
S^{n-2}\times\Sigma_\gamma
\]
admits metrics in \(\mathcal M^+\), but no critical metric for
\(\mathcal R_2\).
\end{enumerate}
\end{corollary}

\begin{proof}
Let \(g\in\mathcal M^+(M)\) be critical for \(\mathcal R_2\).
By Theorem~\ref{Thm:principal} and
Corollary~\ref{cor:virtual-rank-intro},
\(\operatorname{Ric}_{g_u}\geq0\), and \(\pi_1(M)\) is virtually
abelian of rank \(q\leq n-2\).

\smallskip
\noindent\textup{(i)}
By the Cheeger--Gromoll structure theorem
\cite{CheegerGromoll1971}, the universal Riemannian cover splits as
\[
(\widetilde M,\widetilde g_u)
\cong (N^{n-q},g_N)\times\mathbb R^q,
\]
where \(N\) is closed and simply connected.
Since \(n-q\geq2\) and
\(H_{n-q}(N;\mathbb Z_2)\cong\mathbb Z_2\),
the manifold \(N\) is not contractible.
As \(\widetilde M\) is homotopy equivalent to \(N\),
it is not contractible either. Thus \(M\) is not aspherical.

\smallskip
\noindent\textup{(ii)}
If \(\pi_1(M)\) is infinite, then \(q\geq1\).
Choose a connected finite covering \(p:\widehat M\to M\)
with \(\pi_1(\widehat M)\cong\mathbb Z^q\).
As shown in the proof of
Corollary~\ref{cor:virtual-rank-intro},
the metric \(p^*g_u\) admits \(q\) linearly independent
parallel one-forms.
In particular, \(\widehat M\) admits a nowhere vanishing
vector field, so the Poincar\'e--Hopf theorem implies
\[
0=\chi(\widehat M)=\deg(p)\,\chi(M).
\]
Therefore by contraposition, if \(\chi(M)\neq0\) then \(\pi_1(M)\) is finite.

\smallskip
\noindent\textup{(iii)}
Let \(g_{\Sigma_\gamma}\) have constant Gaussian curvature \(-1\).
The product metric
\[
g_K=g_{\mathbb S_K^{n-2}}+g_{\Sigma_\gamma}
\]
has constant scalar curvature
\[
S_{g_K}=(n-2)(n-3)K-2>0
\]
when \(K>2/((n-2)(n-3))\).
Hence \(g_K\in\mathcal M^+(\mathbb S^{n-2}\times\Sigma_\gamma)\).

On the other hand,
\[
\pi_1(\mathbb S^{n-2}\times\Sigma_\gamma)
\cong\pi_1(\Sigma_\gamma)
\]
is not virtually abelian: every connected finite covering
of \(\Sigma_\gamma\) is a closed orientable surface of genus
at least two and therefore has nonabelian fundamental group.
Corollary~\ref{cor:virtual-rank-intro} excludes the existence
of a critical metric for \(\mathcal R_2\) on this product.
\end{proof}

\begin{remark}[Comparison with the classical Laplace eigenvalue problem]
It is useful to compare the preceding conclusions with the classical
picture recalled above. For a volume-normalized Laplace eigenvalue
\[
\overline{\lambda}_k(M,\Delta_g)
=
\lambda_k(M,\Delta_g)\operatorname{Vol}(M,g)^{2/n},
\]
the Euler--Lagrange equation obtained when the whole space of metrics is
varied takes the form, after normalization, 
\[
\sum_{\alpha=1}^{m}
d\varphi_\alpha\otimes d\varphi_\alpha=g,
\]
where $\varphi=(\varphi_\alpha)_\alpha:M\to \mathbb S^{m-1}$ is a minimal homothetic immersion by eigenfunctions of $\Delta_g$ \cite{Nadirashvili1996Berger,ElSoufiIlias2008}. 
In particular, this equation detects every tangent direction and forces the
multiplicity of the relevant eigenvalue to be at least \(n+1\).
In Equation~\eqref{eq:Einstein-harmonic} of
Theorem~\ref{Thm:principal}, the map \(\psi\) need not be an immersion, since
\(d\psi\) vanishes precisely on \(\ker\operatorname{Ric}_{g_u}\).
The multiplicity bound therefore depends on the maximal rank
\(r=\max_M\operatorname{rank}\operatorname{Ric}_{g_u}\),
giving \(\dim E_{\lambda_2}(g)\ge r+1\) instead of \(n+1\).
In particular, \(\psi\) is constant along any connected
Ricci-flat factor of a product metric.
This is already visible on the equality case of
Corollary~\ref{cor:rank-r-intro}.
Let \(2\leq r<n\) and let
\(g=g_{\mathbb S_K^r}+h_F\) on
\(\mathbb S_K^r\times F^{n-r}\), where
\((F,h_F)\) is a closed connected Ricci-flat manifold satisfying
\[
\lambda_1(F,\Delta_{h_F})>rK.
\]
The projection onto the spherical factor, followed by rescaling
to the unit sphere,
\(\psi:\mathbb S_K^r\times F\to\mathbb S^r\), satisfies
\[
\psi^*g_{\mathbb S^r}
=
K g_{\mathbb S_K^r}
=
\frac{1}{r-1}\operatorname{Ric}_g.
\]
Thus \(g\) is critical for \(\mathcal R_2\) by
Theorem~\ref{Thm:principal}.
By contrast, when \(r<n\), \(g\) cannot be critical for the first
volume-normalized positive Laplace eigenvalue. Indeed, since
\(\lambda_1(F,\Delta_{h_F})>rK=\lambda_1(\mathbb S_K^r)\),
all first eigenfunctions are constant along \(F\), so their
differentials cannot recover the metric in those directions.

Thus the harmonic--Einstein equation \eqref{eq:Einstein-harmonic} is less rigid from the extrinsic point of view, since it allows loss of rank and invisible Ricci-flat directions. At the same time, it is more restrictive intrinsically, since its right-hand side is nonnegative and therefore
forces \(\operatorname{Ric}_{g_u}\ge0\), leading to the topological obstructions above. No corresponding nonnegative Ricci curvature conclusion follows from a minimal isometric immersion by Laplace
eigenfunctions.
\end{remark}

\section{Critical metrics for the GJMS ratio of eigenvalues}\label{sec:GJMS}

We now use the tensorial characterization obtained in Proposition~\ref{prop:Euler-Lagrange-GJMS} to derive the higher-order analogue of the harmonic--Einstein system. For this purpose, we recall the curvature tensors associated with \(Q^{(s)}\), following the variational framework of Gover and {\O}rsted \cite{GoverOrsted2013} and Case, Lin and Yuan \cite{CaseLinYuan2019}, and introduce the reduced energy associated with \(P_g^s\).

We assume that \(2s<n\) and that
the conditions~\eqref{conditionsGJMS} hold.
\subsection{Higher order $J$ tensors}
The GJMS operator of order $2s$ satisfies
the conformal covariance law
\begin{equation}\label{eq:GJMS-covariance}
P^s_{u^{4/(n-2s)}g}(\varphi)
=
u^{-\frac{n+2s}{n-2s}}P_g^s(u\varphi),
\end{equation}
and its associated $Q^{(s)}$-curvature is defined by
\begin{equation}\label{eq:Qs-definition}
P_g^s(1)=\frac{n-2s}{2}Q_g^{(s)}.
\end{equation}
Let
\[
\Gamma_g^{(s)}[h]:=D_gQ^{(s)}[h]
\]
be the linearization of $Q^{(s)}$ at $g$, and denote by
$\bigl(\Gamma_g^{(s)}\bigr)^*$ its $L^2$-formal adjoint.
We define the higher $Q$-Ricci tensor by
\begin{equation}\label{eq:Js-definition}
J_g^{(s)}
:=-\frac1s\bigl(\Gamma_g^{(s)}\bigr)^*(1).
\end{equation}
We have
\begin{equation}\label{eq:Js-identities}
\operatorname{tr}_gJ_g^{(s)}=Q_g^{(s)},
\qquad
\operatorname{div}_gJ_g^{(s)}=\frac1{2s}dQ_g^{(s)},
\end{equation}
and therefore the higher order Einstein tensor
\begin{equation}\label{eq:Qs-Einstein}
E_g^{Q^{(s)}}:=-sJ_g^{(s)}+\frac12Q_g^{(s)} g
\end{equation}
is divergence-free. 

The trace and divergence identities follow from the general formulas of Case, Lin and Yuan \cite[Section 3]{CaseLinYuan2019}. For \(s=2\), \(J_g^{(2)}\)
coincides with the \(J\)-tensor of Lin and Yuan
\cite{LinYuanJTensor2017}. See also \cite{GoverOrsted2013,LinYuan2016,LinYuan2022} for more details.

With this normalization, the first variation of the total \(Q^{(s)}\)-curvature takes the form
\begin{align}\label{eq:GJMS-total-Q-variation}
D_g\left(\int_M Q_g^{(s)}\,dv_g\right)[h]
&=
\int_M
\left\langle
\bigl(\Gamma_g^{(s)}\bigr)^*(1)
+\frac12Q_g^{(s)} g,h
\right\rangle_g\,dv_g
\notag\\
&=
\int_M\langle E_g^{Q^{(s)}},h\rangle_g\,dv_g.
\end{align}

As in the conformal Laplacian case, we introduce the reduced GJMS operator defined as

\begin{equation}\label{eq:reduced-GJMS}
L_g^{(s),0}:=P_g^s-P_g^s(1),
\end{equation}
where $P_g^s(1)$ acts by multiplication.
For a smooth map
$\psi=(\psi_1,\ldots,\psi_m):M\to\mathbb S^{m-1}$,
let
\begin{equation}\label{eq:reduced-GJMS-energy}
\mathcal E^{(2s),0}(\psi,g)
:=
\sum_{\alpha=1}^m\int_M
\psi_\alpha L_g^{(s),0}\psi_\alpha\,dv_g.
\end{equation}
Since $|\psi|^2=1$, this can also be written as
\[
\mathcal E^{(2s),0}(\psi,g)
=
\sum_{\alpha=1}^m
\int_M\psi_\alpha P_g^s\psi_\alpha\,dv_g
-\int_MP_g^s(1)\,dv_g.
\]
We denote by $H_g^{(2s)}(\psi)$ its stress tensor
with respect to the metric, that is
\begin{equation}\label{eq:GJMS-stress-definition}
D_g\mathcal E^{(2s),0}(\psi,g)[h]
=
\int_M\langle H_g^{(2s)}(\psi),h\rangle_g\,dv_g,
\end{equation}
where $\psi$ is kept fixed when the metric varies.

Finally, we recall from Section \ref{sec:Clarke} that
\[
\mathcal Q_g^{(s)}(f)=\int_M fP_g^sf\,dv_g,
\qquad
\mathcal N_g(f)=\int_M f^2\,dv_g,
\]
and that \(E_{\lambda_i}(g)\) denotes the eigenspace of \(P_g^s\)
associated with \(\lambda_i(M,P_g^s)\).

Let
\[
\mathcal R_{2s}(g)
:=
\frac{\lambda_2(M,P_g^s)}{\lambda_1(M,P_g^s)},
\]
and let $u$ be a positive first eigenfunction of $P^s_g$.
\subsection{Proof of Theorem \ref{thm:GJMS-ratio}}
Using Euler--Lagrange equations given in Proposition~\ref{prop:Euler-Lagrange-GJMS},
we prove Theorem \ref{thm:GJMS-ratio}, which we restate for convenience.

\begin{theorem*}\label{thm:GJMS-ratio-body}
Let $2s<n$ and assume that the conditions
\eqref{conditionsGJMS} are satisfied.
Then $g$ is Clarke critical for $\mathcal R_{2s}$
if and only if there exist second eigenfunctions
$v_1,\ldots,v_m$ such that
\[
\psi=
\left(\frac{v_1}{u},\ldots,\frac{v_m}{u}\right)
:M\longrightarrow\mathbb S^{m-1}
\]
and, for
\[
g_u=u^{4/(n-2s)}g,
\]
one has
\begin{equation}\label{eq:GJMS-map-equation}
L_{g_u}^{(s),0}\psi
=
\bigl(\mathcal R_{2s}(g)-1\bigr)
\frac{n-2s}{2}Q_{g_u}^{(s)}\psi,
\end{equation}
together with the tensorial equation
\begin{equation}\label{eq:GJMS-metric-equation}
H_{g_u}^{(2s)}(\psi)
=
\bigl(\mathcal R_{2s}(g)-1\bigr)
\frac{n-2s}{2}E_{g_u}^{Q^{(s)}}.
\end{equation}
\end{theorem*}

\begin{proof}[Proof of Theorem \ref{thm:GJMS-ratio}]
Write $\lambda_i=\lambda_i(M,P_g^s)$.
Assume first that $g$ is Clarke critical for $\mathcal R_{2s}$.
By Proposition~\ref{prop:Euler-Lagrange-GJMS}, we may normalize
$u$ by $\int_Mu^2\,dv_g=1$ and choose second eigenfunctions
$v_1,\ldots,v_m$ satisfying~\eqref{eq:GJMS-Clarke-balance}.
Taking the trace and using~\eqref{eq:GJMS-stress-trace}, we obtain
\[
-s\sum_{\alpha=1}^m v_\alpha^2=-su^2,
\]
and therefore
\begin{equation}\label{eq:GJMS-sphere-identity}
\sum_{\alpha=1}^m v_\alpha^2=u^2.
\end{equation}
Hence the map
\[
\psi=
\left(\frac{v_1}{u},\ldots,\frac{v_m}{u}\right)
\]
takes values in $\mathbb S^{m-1}$.
Set
$
g_u=u^{4/(n-2s)}g.
$
Applying~\eqref{eq:GJMS-covariance} first to $1$
and then to $\psi_\alpha=v_\alpha/u$, we get
\begin{equation}\label{eq:GJMS-first-gauge}
P_{g_u}^s(1)
=
u^{-\frac{n+2s}{n-2s}}P_g^s(u)
=
\lambda_1u^{-\frac{4s}{n-2s}},
\end{equation}
and
\begin{equation}\label{eq:GJMS-second-gauge}
P_{g_u}^s\psi_\alpha
=
u^{-\frac{n+2s}{n-2s}}P_g^s(v_\alpha)
=
\lambda_2u^{-\frac{4s}{n-2s}}\psi_\alpha.
\end{equation}
Hence
\[
P_{g_u}^s\psi
=
\frac{\lambda_2}{\lambda_1}P_{g_u}^s(1)\psi.
\]
Subtracting $P_{g_u}^s(1)\psi$ and using
$P_{g_u}^s(1)=\frac{n-2s}{2}Q_{g_u}^{(s)}$
implies
\[
L_{g_u}^{(s),0}\psi
=
\bigl(\mathcal R_{2s}(g)-1\bigr)
\frac{n-2s}{2}Q_{g_u}^{(s)}\psi,
\]
which is~\eqref{eq:GJMS-map-equation}.
Since $|\psi|^2=1$, the right-hand side
is normal to the sphere, and the tangential part
of $L_{g_u}^{(s),0}\psi$ vanishes.
Thus $\psi$ is a critical point of the reduced
GJMS energy under the constraint $|\psi|=1$.

It remains to identify the tensorial equation.
For this, we first replace ~\eqref{eq:GJMS-stress-decomposition}
in~\eqref{eq:GJMS-Clarke-balance} to obtain,
\[
\sum_{\alpha=1}^mD_g\mathcal N_g(v_\alpha)[h]
=
\frac12\int_M
\left(\sum_{\alpha=1}^m v_\alpha^2\right)
\operatorname{tr}_gh\,dv_g
=
D_g\mathcal N_g(u)[h].
\]
Therefore for every metric variation $h$ we have,
\begin{equation}\label{eq:GJMS-quadratic-balance-g}
\sum_{\alpha=1}^m
D_g\mathcal Q_g^{(s)}(v_\alpha)[h]
=
\mathcal R_{2s}(g)D_g\mathcal Q_g^{(s)}(u)[h].
\end{equation}
Using the conformal covariance 
\begin{equation}\label{eq:GJMS-quadratic-form-covariance}
\mathcal Q_g^{(s)}(u\varphi)
=
\mathcal Q_{g_u}^{(s)}(\varphi),
\end{equation}
and differentiating~\eqref{eq:GJMS-quadratic-form-covariance}
in~\eqref{eq:GJMS-quadratic-balance-g}
gives, for every $ h_u=u^{4/(n-2s)}h$ with $u$ fixed,
\begin{equation}\label{eq:GJMS-quadratic-balance-gu}
\sum_{\alpha=1}^m
D_{g_u}\mathcal Q_{g_u}^{(s)}(\psi_\alpha)[h_u]
=
\mathcal R_{2s}(g)
D_{g_u}\mathcal Q_{g_u}^{(s)}(1)[h_u].
\end{equation}
Since $|\psi|^2=1$, this is exactly
\begin{equation}\label{eq:GJMS-reduced-energy-variation}
D_{g_u}\mathcal E^{(2s),0}(\psi,g_u)[h_u]
=
\bigl(\mathcal R_{2s}(g)-1\bigr)
D_{g_u}\mathcal Q_{g_u}^{(s)}(1)[h_u].
\end{equation}

Recalling that
$\mathcal Q_g^{(s)}(f)=\int_M fP_g^sf\,dv_g$,
we obtain from~\eqref{eq:Qs-definition}
\begin{equation}\label{eq:GJMS-Q-total}
\mathcal Q_{g_u}^{(s)}(1)
=
\int_M P_{g_u}^s(1)\,dv_{g_u}
=
\frac{n-2s}{2}
\int_M Q_{g_u}^{(s)}\,dv_{g_u}.
\end{equation}

Combining~\eqref{eq:GJMS-reduced-energy-variation}
with~\eqref{eq:GJMS-Q-total}
and~\eqref{eq:GJMS-total-Q-variation},
and using the definition~\eqref{eq:GJMS-stress-definition},
we obtain
\[
H_{g_u}^{(2s)}(\psi)
=
\bigl(\mathcal R_{2s}(g)-1\bigr)
\frac{n-2s}{2}E_{g_u}^{Q^{(s)}},
\]
which is~\eqref{eq:GJMS-metric-equation}.

Finally, from~\eqref{eq:GJMS-metric-equation}
and
$\operatorname{tr}_{g_u}E_{g_u}^{Q^{(s)}}
=\frac{n-2s}{2}Q_{g_u}^{(s)}$, we get
\[
\operatorname{tr}_{g_u}H_{g_u}^{(2s)}(\psi)
=
\bigl(\mathcal R_{2s}(g)-1\bigr)
\frac{(n-2s)^2}{4}Q_{g_u}^{(s)}.
\]
Subtracting $1/(n-2s)$ times this trace
from the tensorial equation and multiplying by $-1/s$ we finally deduce the  trace-reversed form
\begin{equation}\label{eq:GJMS-trace-reversed}
-\frac1s\left(
H_{g_u}^{(2s)}(\psi)
-\frac{\operatorname{tr}_{g_u}H_{g_u}^{(2s)}(\psi)}
{n-2s}g_u
\right)
=
\bigl(\mathcal R_{2s}(g)-1\bigr)
\frac{n-2s}{2}J_{g_u}^{(s)}.
\end{equation}

Conversely, assume that the eigenfunctions
$u,v_1,\ldots,v_m$ exist, that $\psi$ is sphere-valued,
and that~\eqref{eq:GJMS-metric-equation} holds.
Equation~\eqref{eq:GJMS-total-Q-variation}
allows us to reverse the preceding argument:
the tensorial equation gives
\eqref{eq:GJMS-reduced-energy-variation},
hence~\eqref{eq:GJMS-quadratic-balance-gu},
and then, by conformal covariance,
\eqref{eq:GJMS-quadratic-balance-g}.
Since $\sum_\alpha v_\alpha^2=u^2$,
the norm-variation terms
in~\eqref{eq:GJMS-stress-decomposition}
cancel and we recover~\eqref{eq:GJMS-Clarke-balance}.
Proposition~\ref{prop:Euler-Lagrange-GJMS}
now implies that $g$ is Clarke critical
for $\mathcal R_{2s}$.
\end{proof}
\begin{remark}

For $s=2$, corresponding to the Paneitz operator, the two equations become
\begin{equation}\label{eq:Paneitz-system}
L_{g_u}^{(2),0}\psi
=
\bigl(\mathcal R_4(g)-1\bigr)
\frac{n-4}{2}Q_{g_u}^{(2)}\psi,
\qquad
H_{g_u}^{(4)}(\psi)
=
\bigl(\mathcal R_4(g)-1\bigr)
\frac{n-4}{2}E_{g_u}^{Q^{(2)}}.
\end{equation}
where $g_u=u^{4/(n-4)}g$.
The trace-reversed equation reduces to
\[
-\frac12\left(
H_{g_u}^{(4)}(\psi)
-\frac{\operatorname{tr}_{g_u}H_{g_u}^{(4)}(\psi)}
{n-4}g_u
\right)
=
\bigl(\mathcal R_4(g)-1\bigr)
\frac{n-4}{2}J_{g_u}^{(2)}.
\]
This identity highlights the role of $J_{g_u}^{(2)}$
as the higher-order counterpart of the Ricci tensor
in~\eqref{eq:Einstein-harmonic}. 

This is the fourth-order analogue of the harmonic map equation into a sphere. Since $|\psi|^2=1$, the right-hand side is normal to the sphere, so the tangential part of $L_{g_u}^{(2),0}\psi$ vanishes:
$$
\Pi_{T_\psi\mathbb S^{m-1}}
\left(
L_{g_u}^{(2),0}\psi
\right)=0.
$$
Thus $\psi$ is a critical point of the reduced Paneitz energy under the constraint $|\psi|=1$. This is the Paneitz analogue of harmonicity. It is not the usual biharmonic map equation, the operator $L_{g_u}^{(2),0}$ has leading term $\Delta_{g_u}^2$, but also contains the lower-order curvature terms required by conformal covariance. The second equation in \eqref{eq:Paneitz-system} is the fourth-order analogue of the harmonic--Einstein equation \eqref{eq:Einstein-harmonic}.
\end{remark}

\bigskip

\subsection*{Acknowledgements}

The author is grateful to David Tewodrose for his support, his valuable suggestions and his careful reading of the manuscript. This work is funded by the Research Foundation– Flanders (FWO) via the Odysseus II programme no. G0DBZ23N.

\bibliographystyle{plain}
\bibliography{harmonic_einstein_references}
\end{document}